\documentclass[12pt]{amsart}
\usepackage[all, cmtip]{xy}

\usepackage{times,amsfonts,amsmath,amstext,amsbsy,amssymb,
amsopn,amsthm,upref,eucal}

\usepackage{amsmath,amssymb,amscd}

\usepackage{amssymb, amsmath, amscd, amsthm, amsfonts, amscd, mathrsfs}

\usepackage{url}
\newcommand \spann {{\mathrm{span}}}
\newcommand \Mat {{\mathrm{Mat}}}
\newcommand \Prob {{\mathbb P}}

\newtheorem{theorem}{Theorem}

\newtheorem*{teorema}{Theorem}

\newtheorem{assumption}{Assumption}
\newtheorem{corollary}{Corollary}
\newtheorem{proposition}{Proposition}

\DeclareMathOperator{\Conf}{Conf}
\DeclareMathOperator{\tr}{tr}

\newcommand{\conf}{\mathrm{conf}}

\newcommand\scrI{{\mathscr I}}

\begin{document}

\title
{Infinite Determinantal Measures}
\author{Alexander I. Bufetov}

\address{Laboratoire d'Analyse, Topologie, Probabilit\'es, CNRS, Marseille}
\address{The Steklov Institute of Mathematics, Moscow}
\address {The Institute for Information Transmission Problems, Moscow}
\address{National Research University Higher School of Economics, Moscow}
\address{The Independent University of Moscow}
\address{Rice University, Houston}
\date{}
\maketitle

\begin{abstract}
Infinite determinantal measures introduced in this note are
inductive limits of determinantal measures on an exhausting
family of subsets of the phase space. Alternatively,
an infinite determinantal measure can be described
as a product of a determinantal point process
and a convergent, but not integrable, multiplicative functional.

Theorem 2, the main result announced in this note, gives an explicit description for
the ergodic decomposition of infinite Pickrell measures on the
spaces of infinite complex matrices in terms of  infinite
determinantal measures obtained by finite-rank perturbations of
Bessel point processes.

\end{abstract}

\section{Introduction}
\subsection{Outline of the main results.}
In this section, our aim is to construct sigma-finite analogues of determinantal measures on spaces of configurations. In Theorem \ref{mainthm} of Section 4, infinite determinantal measures will be seen to arise  in the ergodic decomposition of infinite 
unitarily-invariant measures on spaces of infinite complex matrices. 

Informally, a configuration on the phase space $E$ is an unordered collection of points (called {\it particles}) of $E$, possibly with multiplicities; the main assumption is that a bounded subset of $E$ contain only finitely many particles of a given configuration.

To a function $g$ on $E$ assign its {\it multiplicative functional} $\Psi_g$ on the space of configurations: the functional $\Psi_g$ is obtained by multiplying the values of $g$ over all particles of a configuration (see (\ref{mult-fun-def})). 
 A probability measure on the space of 
configurations on $E$ is uniquely characterized
by prescribing the expectations of multiplicative functionals; for {\it determinantal} probability measures these expectations are given by special Fredholm determinants, see e.g. \cite{Soshnikov}; 
 the definition is also recalled in  
(\ref{eq1}) below. 

Given a subset $E^{\prime}\subset E$, consider the subset 
${\mathrm {Conf}}(E, E^{\prime})$ of those configurations whose all  particles lie in $E^{\prime}$; in Proposition \ref{indsubset} below, we shall see that under some additional asumptions the restriction of a determinantal point process onto ${\mathrm {Conf}}(E, E^{\prime})$ is again determinantal.

Our main example, the measure ${\mathbb B}^{(s)}$ of (\ref{bs-def}), is defined on the space of configurations on $(0, +\infty)$. Almost every configuration is infinite and bounded according to ${\mathbb B}^{(s)}$;
the particles  accumulate at zero. 
If one  takes $R>0$ and requires all particles to lie in $(0,R)$,
 then the  induced measure of ${\mathbb B}^{(s)}$ on the resulting subset of configurations 
is finite, and, after normalization, determinantal. As $R$ goes to infinity, the measure of the subset 
${\mathrm {Conf}}((0, +\infty); (0, R))$  grows, and the measure of the space of all configurations 
is infinite.

Our general construction will similarly exhaust $E$ by subsets $E_n$ in such a way that the weight of ${\mathrm {Conf}}(E; E_n)$ is positive and finite, and the normalized restriction of our infinite determinantal measure onto the subset  ${\mathrm {Conf}}(E; E_n)$ is determinantal. A simple example is given by ``infinite orthogonal polynomial ensembles'', see  (\ref{ore}) below. The measure  ${\mathbb B}^{(s)}$ is a scaling limit of such ensembles. 
We proceed to precise formulations.

\subsection{Construction of infinite determinantal measures}
Let $E$ be a locally compact complete metric space, and let $\Conf(E)$ be the space of configurations on $E$ endowed with the natural Borel
structure (see, e.g., \cite {Lenard},~\cite{Soshnikov}).

Given a Borel subset $E'\subset E$, we let $\Conf(E, E')$ be the subspace of
configurations all whose particles lie in $E'$.

Given a measure $\mathbb{B}$ on a set $X$ and a measurable subset $Y\subset X$ such
that $0<\mathbb{B}(Y)<+\infty$, we let $\mathbb{B}\left|_Y\right.$ stand for the
restriction of the measure $\mathbb{B}$ onto the subset $Y$.

An {\it infinite determinantal measure} is a  $\sigma$-finite Borel measure $\mathbb{B}$
on $\Conf(E)$ admitting a filtration of the space $E$ by Borel subsets $E_n$,
$n\in\mathbb{N}$:

$$E_1\subset E_2\subset \ldots\subset E_n\subset \ldots\;,\;\;\bigcup\limits_{n=1}^{\infty}E_n=E$$
such that for any $n\in\mathbb{N}$ we have
\begin{enumerate}
\item $0<\mathbb{B}\left(\Conf(E, E_n)\right)<+\infty$;
\item the normalized restriction
$$\frac{\mathbb{B}\left|_{\Conf(E, E_n)}\right. }{\mathbb{B}\left(\Conf(E, E_n)\right)}$$
is a determinantal measure;
\item $\mathbb{B}\left(\Conf(E)\backslash\bigcup\limits_{n=1}^{\infty}(\Conf(E,
E_n)\right)=0\,.$
\end{enumerate}
%One can construct infinite determinantal measures in the following way.

Let $\mu$ be a $\sigma$-finite Borel measure on $E$. 
By the Macch{\`i}-Soshnikov Theorem, under some additional assumptions, a determinantal measure can be assigned to an operator of orthogonal projection, or, in other words,  to a closed subspace of $L_2(E, \mu)$. In a similar way, an infinite determinantal measure will be assigned to a subspace $H$  of {\it locally} square-integrable functions. For example, for  infinite analogues of orthogonal polynomial ensembles, $H$ is the subspace of weighted polynomials, see Subsection 1.3 below.

Let $L_{2,\mathrm{loc}}(E,\mu)$ be the
space of measurable functions on $E$, locally square integrable with respect to
$\mu$, let $\scrI_1(E,\mu)$ be the space of trace-class operators in $L_2(E,\mu)$
and let $\scrI_{1,\mathrm{loc}}(E,\mu)$ be the space of operators
on $L_2(E,\mu)$ that are locally of trace class (precise definitions are recalled in Section 2).

Let $H\subset L_{2,\mathrm{loc}}(E,\mu)$ be a linear subspace. If
$E'\subset E$ is a Borel subset such that $\chi_{E'}H$ is a closed subspace of
$L_2(E,\mu)$, then we denote by $\Pi^{E'}$ the operator of
orthogonal projection onto the subspace $\chi_{E'}H\subset
L_2(E,\mu).$ We now fix a Borel subset $E_0\subset E$;  informally, $E_0$ is the set where the particles accumulate. We impose
the following assumption on $E_0$ and $H$.
\begin{assumption}\label{he}
\begin{enumerate}
\item For any bounded Borel set $B\subset E$, the space  $\chi_{E_0\cup B} H$ is a closed subspace of $L_2(E,\mu)$;
 \item For any bounded Borel set $B\subset E\setminus E_0$, we have
\begin{equation}
\Pi^{E_0\cup B}\in \scrI_{1,\mathrm{loc}}(E,\mu),\quad \chi_B\Pi^{E_0\cup B}\chi_B\in\scrI_{1}(E,\mu);
\end{equation}
\item If $\varphi\in H$ satisfies $\chi_{E_0}\varphi =0$, then $\varphi=0.$
\end{enumerate}
\end{assumption}

\begin{theorem}\label{infdet-he}
Let $E$ be a locally compact complete metric space, and let $\mu$ be
a $\sigma$-finite Borel measure on $E$. If a subspace
$H\subset L_{2,\mathrm{loc}}(E,\mu)$ and a Borel subset
$E_0\subset E$ satisfy Assumption \ref{he}, then there exists a
$\sigma$-finite Borel measure $\mathbb{B}$ on $\mathrm{\Conf}(E)$
such that
\begin{enumerate}
\item ${\mathbb B}$-almost every configuration has at most finitely
many particles outside of $E_0$;
\item for any bounded Borel (possibly empty) subset
$B\subset E\setminus E_0$ we have $0<\mathbb{B}(\mathrm{\Conf}(E;E_0\cup B))<+\infty$
and  $$\frac{\mathbb{B}|_{\mathrm{\Conf}(E;E_0\cup
B)}}{\mathbb{B}({\mathrm{\Conf}(E;E_0\cup
B)})}=\mathbb{P}_{\Pi^{E_0\cup B}}.$$
\end{enumerate}
The requirements (1) and (2) determine the measure $\mathbb{B}$
uniquely up to multiplication by a positive constant.
\end{theorem}
We denote ${\mathbf B}(H,E_0)$ the one-dimensional cone of nonzero infinite determinantal measures induced by
$H$ and $E_0$, and, slightly abusing notation, we write ${\mathbb B}={\mathbb B}(H, E_0)$
for a representative of the cone.

{\bf {Remark.}} If $B$ is a bounded set, then, by definition, we have
$${\mathbf B}(H,E_0)={\mathbf B}(H,E_0\cup B).$$

{\bf {Remark.}}  If $E^{\prime}\subset E$ is a Borel subset such that $\chi_{E_0\cup E^{\prime}}$ is
a closed subspace in $L_2(E, \mu)$ and the operator $\Pi^{E_0\cup E^{\prime}}$ of orthogonal projection onto the subspace
$\chi_{E_0\cup E^{\prime}}H$ satisfies
\begin{equation}
\Pi^{E_0\cup E^{\prime}}\in \scrI_{1,\mathrm{loc}}(E,\mu),\quad \chi_{E^{\prime}}\Pi^{E_0\cup E^{\prime}}\chi_{E^{\prime}}\in\scrI_{1}(E,\mu),
\end{equation}
then, exhausting $E^{\prime}$ by bounded sets, from Theorem \ref{infdet-he} one easily obtains
$0<\mathbb{B}(\mathrm{\Conf}(E;E_0\cup E^{\prime}))<+\infty$
and
$$
\frac{\mathbb{B}|_{\mathrm{\Conf}(E;E_0\cup
E^{\prime})}}{\mathbb{B}({\mathrm{\Conf}(E;E_0\cup
E^{\prime})})}=\mathbb{P}_{\Pi^{E_0\cup E^{\prime}}}.
$$
\subsection{Infinite orthogonal polynomial ensembles}
Take an interval $[a,b)$ in ${\mathbb R}$, let $\mathrm{Leb}=dx$ on $[a,b)$ be the Lebesgue measure on $[a,b)$, let $\rho$ be a positive continuous function on $[a,b)$, and assume
$\int_a^b\rho(x)dx=+\infty$.
 Take $n\in\mathbb{N}$ and
endow the set $[a,b]^N$ with  the measure 
\begin{equation}\label{ore}
\prod\limits_{1\leqslant i,j\leqslant N}(x_i-x_j)^2\prod
\limits_{i=1}^N\rho(x_i)dx_i,
\end{equation}
an infinite analogue of an orthogonal polynomial ensemble.

 For any $b_1\in[a,b)$, the induced
measure
\begin{equation}\label{ore-ind}
\prod\limits_{1\leqslant i<
j\leqslant N}(x_i-x_j)^2\prod\limits_{i=1}^N\rho(x_i)\chi_{[a,b_1]}
(x_i)dx_i
\end{equation}
is finite and, after normalization, can be represented in
determinantal form
$$\frac{1}{N!}\det K_N^{\rho,b_1}(x_i,x_j)\prod\limits_{i=1}^N
\rho(x_i)\chi_{[a,b_1]}(x_i)dx_i,$$
where $K_N^{\rho,b_1}$ is the $N$-th Christoffel-Darboux kernel formed by orthonormal polynomials
corresponding to the ``induced'' weight $\rho(x)\chi_{[a,b_1]}(x).$

The infinite measure (\ref{ore}) is thus an
infinite
determinantal measure corresponding to the subspace $H\subset L_{2, \mathrm{loc}}([a,b), \mathrm{Leb})$ spanned by the functions $x^k\sqrt{\rho(x)}$, $k=0, \dots, N-1$, and the subset $E_0=[a, b_1)$ for an arbitrary $b_1\in (a,b)$. In the problem of ergodic decomposition of infinite Pickrell measures we shall be especially  interested in studying scaling limits of such ``infinite orthogonal polynomial ensembles''.
\subsection{Organization of the paper.}
%The remainder of this note is organized as follows. 
In the next subsection it is shown 
that, under certain additional assumptions, an infinite determinantal 
measure times a multiplicative functional yields after normalization  a determinantal point process; 
for determinantal probability measures this has been established in \cite{Buf-umn}. We then proceed to our main example of infinite determinantal measures, namely, those obtained as finite-rank perturbations of 
determinantal point processes. The ergodic decomposition measures of infinite Pickrell measures will be seen to be of this type. In the following subsection it is established that induced processes of an infinite determinantal measure obtained by finite rank perturbation, converge to the unperturbed process. 

In Section 2  we  recall the 
definition of determinantal point processes,  study the properties of multiplicative functionals of these processes, thus extending the results of \cite{Buf-umn}, and give a sketch of the proof of Theorem \ref{infdet-he}.

  %Pickrell measures on the space of infinite matrices are introduced 
  In Section 3 we recall the construction, due to Pickrell \cite{Pickrell1}, \cite{Pickrell2}, \cite{Pickrell3} in the finite case 
  (see also Neretin \cite{Neretin}) and to Borodin and Olshanski \cite{BorOlsh} in the infinite case, of Pickrell measures on the space of infinite matrices. We then recall the Olshanski-Vershik approach (see  \cite{Vershik}, \cite{OlshVershik}) to the  Pickrell classification of finite ergodic unitarily-invariant measures on spaces of infinite matrices as well as the result of \cite{Buf-inferg} that implies that the ergodic components of infinite Pickrell measures are almost surely finite; only the decomposing measure is infinite. 
 
In Section 4 we
 %formulate the main result of the paper, Theorem \ref{mainthm}. We 
 start by considering finite Pickrell mesures, for which the ergodic decomposition is given, up to a change of variable,  by the Bessel point process of Tracy and Widom \cite{TracyWidom}. The main result of the paper, Theorem \ref{mainthm} , then says that 
the ergodic decomposition of infinite Pickrell measures is induced by  infinite determinantal measures obtained as an explicitly  given finite-rank perturbation of the Bessel point processes occurring in the ergodic decomposition of finite Pickrell measures.
The scaling limit argument sketched at the end of the section uses precisely the representation, developed in Section 1,  of  infinite determinantal measures as products of finite determinantal measures and  multiplicative functionals. 
\subsection{Multiplicative functionals}

Let $g$ be a non-negative measurable function on $E$, and introduce the
{\it multiplicative functional} $\Psi_g:\Conf(E)\to\mathbb{R}$ by the formula
 \begin{equation} \label{mult-fun-def}
 \Psi_g(X)=\prod\limits_{x\in X}g(x).
 \end{equation}
If the infinite product
$\prod\limits_{x\in X}g(x)$ absolutely converges to $0$ or to $\infty$, then we set, respectively,
$\Psi_g(X)=0$ or $\Psi_g(X)=\infty$. If the product in the right-hand side fails to converge absolutely,
then the multiplicative functional is not defined.

%As we shall see, the key condition for  integrability is that the subspace $\sqrt{g}H$ lie in $L_2(E, \mu)$.
%Note here that if $P$ is an operator of orthogonal projection onto a
%subspace $L\subset L_2(E, \mu)$

 We start with an auxiliary proposition.
\begin{proposition}\label{gh-closed} Let a subspace $H\subset L_{2,\mathrm{loc}}(E,\mu)$ and a Borel subset $E_0\subset E$
satisfy Assumption \ref{he}. Let $g$ be a positive bounded measurable function on $E$ such that
\begin{enumerate}
\item for any bounded subset $B\subset E$ there exists $\varepsilon_0=\varepsilon_0(B)>0$ such that $g(x)>\varepsilon_0$ for
all $x\in E_0\cup B$;
\item we have $\sqrt{g}H\subset L_2(E,\mu)$.
\end{enumerate}
Then $\sqrt{g}H$ is a closed subspace in $L_2(E,\mu)$.
\end{proposition}
Under the assumptions of Proposition~\ref{gh-closed}, let $\Pi^g$ be the operator of orthogonal projection
onto the closed subspace $\sqrt{g}H$.

Our next aim is to give sufficient conditions for integrability of multiplicative
functionals with respect to infinite determinantal measures.
We restrict ourselves to the case when the function $g$ only takes values in $(0,1]$.

\begin{proposition}\label{mult-he}Let a subspace $H\subset L_{2,\mathrm{loc}}(E,\mu)$
and a Borel subset $E_0\subset E$ satisfy Assumption \ref{he}, and let $g\colon E\to(0,1]$ be a
measurable function such that:
\begin{enumerate}
\item for any bounded subset $B\subset E$ there exists $\varepsilon_0=\varepsilon_0(B)>0$ such that $g(x)>\varepsilon_0$ for
all $x\in E_0\cup B$;
\item $\sqrt{g}H\subset L_2(E,\mu)$;
\item $\sqrt{1-g}\chi_{E_0}\Pi^g\chi_{E_0}\sqrt{1-g}\in\scrI_1(E,\mu)$;
\item $\Pi^g\in \scrI_{1,\mathrm{loc}}(E,\mu)$;
\item $\chi_{E\setminus E_0}\Pi^g\chi_{E\setminus E_0}
\in\scrI_1(E,\mu)$.
\end{enumerate}
Then the multiplicative functional $\Psi_g$ is $\mathbb{B}(H,E_0)$-almost surely positive, and we have
\begin{enumerate}
\item $$\Psi_g\in L_1(\Conf(E),\mathbb{B});$$
\item $$\displaystyle\frac{\displaystyle
\Psi_g\mathbb{B}}{\displaystyle
\int\limits_{\Conf(E)}\Psi_g\,d\mathbb{B}}=\mathbb{P}_{\Pi^g}.$$
\end{enumerate}
\end{proposition}

We can therefore write
$
{\mathbb B}=C\cdot \Psi_{1/g}\cdot \mathbb{P}_{\Pi^g},
$
where $C$ is a positive constant. Our infinite determinantal measure is thus represented as a product of a determinantal probability measure and a convergent non-integrable multiplicative functional.

\subsection{Infinite determinantal measures obtained as finite-rank perturbations of  determinantal probability measures.}

We now consider infinite determinantal measures induced by subspaces $H$ obtained by adding a
finite-dimensional subspace $V$ to a closed subspace $L\subset L_2(E, \mu)$.

Let, therefore,  $Q\in\scrI_{1,\mathrm{loc}}(E,\mu)$ be the operator of orthogonal projection onto a closed subspace $L\subset L_2(E, \mu)$, let $V$ be a finite-dimensional subspace of $L_{2, \mathrm{loc}}(E, \mu)$, and set $H=L+V$.  Let  $E_0\subset E$ be a Borel subset.
We shall need the following assumption on $L, V$ and $E_0$.
\begin{assumption}\label{lve}
\begin{enumerate}
\item $\chi_{E\setminus E_0}Q\chi_{E\setminus E_0}\in\scrI_1(E,\mu)$;
\item $\chi_{E_0}V\subset L_2(E,\mu)$;
\item if $\varphi\in V$ satisfies $\chi_{E_0}\varphi\in \chi_{E_0}L$, then $\varphi=0$;
\item 
if $\varphi\in L$ satisfies $\chi_{E_0}\varphi=0$, then $\varphi=0$.
\end{enumerate}
\end{assumption}

\begin{proposition} \label{lvehe} If  $L$, $V$ and  $E_0$ satisfy Assumption \ref{lve} then
the subspace  $H=L+V$ and  $E_0$ satisfy Assumption \ref{he}.
\end{proposition}

In particular, for any bounded Borel subset $B$, 
the subspace  $\chi_{E_0\cup B}L$ is closed, as  one sees by taking $E^{\prime}=E_0\cup B$ in the following clear 
\begin{proposition}\label{closedsubsp}
Let $Q\in\scrI_{1,\mathrm{loc}}(E,\mu)$ be the operator of orthogonal projection
onto a closed subspace $L\in L_2(E,\mu)$. Let $E'\subset E$ be a Borel subset such that
$\chi_{E'}Q\chi_{E'}\in\scrI_1(E,\mu)$ and that for any function $\varphi\in L$, the equality $\chi_{E'}\varphi=0$ implies $\varphi=0$. Then the subspace $\chi_{E'}L$ is closed
in $L_2(E, \mu)$.\end{proposition}

The subspace $H$ and the Borel subset $E_0$ therefore define an infinite determinantal measure $\mathbb{B}=\mathbb{B}(H,E_0)$.
We now adapt the formulation of Proposition~\ref{mult-he} to this particular case.

\begin{proposition}\label{finite-mult} Let $L$, $V$, and $E_0$ satisfy Assumption~\ref{lve}, let $\mathbb{B}$ be the corresponding
infinite determinantal measure, and let $g\colon E\to (0,1]$ be a positive measurable function.
 If
%\begin{enumerate}
$\sqrt{1-g}Q\sqrt{1-g}\in\scrI_1(E,\mu),$
%\item $\chi_{E\setminus E_0}\Pi\chi_{E\setminus E_0}\in\scrI_1(E,\mu),$
%\end{enumerate}
then the multiplicative functional $\Psi_g$ is $\mathbb{B}$-almost surely well-defined and positive.

 If,  additionally, we assume
 \begin{enumerate}
 \item $\sqrt{g}V\subset L_2(E,\mu)$;
 \item for any bounded subset $B\subset E$ there exists $\varepsilon_0=\varepsilon_0(B)>0$ such that $g(x)>\varepsilon_0$ for
all $x\in E_0\cup B$,
\end{enumerate}
then
\begin{enumerate}
\item $$\Psi_g\in L_1(\Conf(E),\mathbb{B});$$
\item $$\displaystyle\frac{\displaystyle
\Psi_g\mathbb{B}}{\displaystyle
\int\limits_{\Conf(E)}\Psi_g\,d\mathbb{B}}=\mathbb{P}_{\Pi^g},$$
\end{enumerate}
where, as before, $\Pi^g$ is the operator of orthogonal projection onto the closed subspace $\sqrt{g}H$.

%\end{enumerate}
\end{proposition}

{\bf Remark.} The subspace $\sqrt{g}H$ is closed by Proposition \ref{gh-closed}.

\subsection{Convergence of approximating kernels.}
Our next aim is to show that, under certain additional assumptions, if a sequence $g_n$ of measurable functions converges to $1$, then the operators $\Pi^{g_n}$ considered in Proposition \ref{finite-mult}
converge to $Q$ in $\scrI_{1, \mathrm{loc}}(E, \mu)$.

%As above, let $Q\in\scrI_{1,\mathrm{loc}}(E,\mu)$ be the operator of orthogonal projection onto a closed subspace $L\subset L_2(E, \mu)$,
%let $V$ be a finite-dimensional subspace of $L_{2, \mathrm{loc}}(E, \mu)$, let $E_0\subset E$ be a Borel subset, and suppose that Assumption  \ref{lve} is verified.

 %As above, if the subspace $\sqrt{g}H=\sqrt{g}L+\sqrt{g}V$ is closed, then we let  $\Pi^g$ be  the operator of orthogonal projection onto  $\sqrt{g}H$, and,
 %for brevity, we write $\Pi^{E_1}=\Pi^{\chi_{E_1}}$.

Given  two closed subspaces $H_1, H_2$ in $L_2(E, \mu)$, let $\alpha(H_1, H_2)$
be the angle between $H_1$ and $H_2$, defined as the infimum of angles between all nonzero vectors
in $H_1$ and $H_2$; recall that if one of the subspaces has finite dimension, then the infimum is achieved.
\begin{proposition} \label{gn-conv}
Let $L$, $V$, and $E_0$ satisfy Assumption~\ref{lve}, and assume additionally that we have $V\cap L_2(E,\mu)=0.$
Let $g_n:E\to (0,1]$ be a sequence of positive measurable functions such that
\begin{enumerate}
\item for all $n\in {\mathbb N}$ we have
$\sqrt{1-g_n}Q\sqrt{1-g_n}\in\scrI_1(E,\mu)$;\item for all $n\in {\mathbb N}$ we have $\sqrt{g_n}V\subset L_2(E, \mu)$;
\item there exists $\alpha_0>0$ such that for all $n$ we have
$$
\alpha(\sqrt{g_n}H, \sqrt{g_n}V)\geq \alpha_0;
$$
\item for any bounded $B\subset E$  we have
$$
\inf\limits_{n\in {\mathbb N}, x\in E_0\cup B} g_n(x)>0;
$$
$$
\lim\limits_{n\to\infty} \sup\limits_{x\in E_0\cup B} \left|g_n(x)-1\right|=0.
$$
\end{enumerate}
Then, as $n\to\infty$, we have
$$
\Pi^{g_n}\to Q\text{ in }\scrI_{1, \mathrm{loc}}(E, \mu).
$$
\end{proposition}

Using the second remark after Theorem \ref{infdet-he}, one can extend Proposition \ref{gn-conv} also to nonnegative functions that admit zero values. Here we restrict ourselves to characteristic functions of the form $\chi_{E_0\cup B}$ with $B$ bounded, in which case we have the following

\begin{corollary}
\label{cor-conv}
Let $B_n$ be an increasing sequence of bounded Borel sets exhausting $E\setminus E_0$.
 If there exists $\alpha_0>0$ such that for all $n$ we have
$$
\alpha(\chi_{E_0\cup B_n}H, \chi_{E_0\cup B_n}V)\geq \alpha_0,
$$
then
$$
\Pi^{E_0\cup B_n}\to Q\text{ in }\scrI_{1,\mathrm{loc}}(E, \mu).
$$
\end{corollary}

Informally, Corollary \ref{cor-conv} means that, as $n$ grows, the
induced processes of our determinantal measure
on subsets $\Conf(E; E_0\cup B_n)$ converge to the ``unperturbed''
determinantal point process $\Prob_{Q}$.

\section{Multiplicative Functionals of Determinantal Point Processes}

\subsection{Locally integrable functions and locally trace class operators}

Recall that $L_{2,  \mathrm{loc}}(E, \mu)$ is the space of all measurable functions $f:E\to {\mathbb C}$
such that for any bounded subset $B\subset E$ we have
\begin{equation}
\label{bintef}
\int\limits_B |f|^2d\mu<+\infty.
\end{equation}

Choosing an exhausting family $B_n$ of bounded sets (for instance,
balls of radius tending to infinity) and using (\ref{bintef}) with
$B=B_n$, we endow the space $L_{2, \mathrm{loc}}(E, \mu)$ with a countable
family of seminorms which turns it into a complete separable metric
space; the topology thus defined does not, of course, depend on the
specific choice of the exhausting family.

Let $\scrI_{1}(E,\mu)$ be the ideal of trace class operators
${\widetilde K}\colon L_2(E,\mu)\to L_2(E,\mu)$ (see volume~1 of~\cite{ReedSimon} for
the precise definition); the symbol
$||{\widetilde K}||_{\scrI_{1}}$ will stand for the
$\scrI_{1}$-norm of the operator ${\widetilde K}$.
Let $\scrI_{2}(E,\mu)$ be the
ideal of Hilbert-Schmidt operators ${\widetilde K}\colon
L_2(E,\mu)\to L_2(E,\mu)$; the symbol $||{\widetilde
K}||_{\scrI_{2}}$ will stand for the $\scrI_{2}$-norm of
the operator ${\widetilde K}$.

Let  $\scrI_{1,  \mathrm{loc}}(E,\mu)$ be the space of operators  $K\colon L_2(E,\mu)\to L_2(E,\mu)$
such that for any bounded Borel subset $B\subset E$
we have $$\chi_BK\chi_B\in\scrI_1(E,\mu).$$
Again, we endow the space $\scrI_{1,  \mathrm{loc}}(E,\mu)$
with a countable family of semi-norms
\begin{equation}
\label{btrcl}
||\chi_BK\chi_B||_{\scrI_1}
\end{equation}
where, as before, $B$ runs through an exhausting family $B_n$ of bounded sets.

\subsection{Determinantal Point Processes}

A Borel probability measure $\mathbb{P}$ on
$\Conf(E)$ is called
\textit{determinantal} if there exists an operator
$K\in\scrI_{1,  \mathrm{loc}}(E,\mu)$ such that for any bounded measurable 
function $g$, for which $g-1$ is supported in a bounded set $B$, 
we have 
\begin{equation}
\label{eq1}
\mathbb{E}_{\mathbb{P}}\Psi_g
=\det\biggl(1+(g-1)K\chi_{B}\biggr).
\end{equation}
The Fredholm determinant in~\eqref{eq1} is well-defined since
$K\in \scrI_{1, \mathrm{loc}}(E,\mu)$.
The equation (\ref{eq1}) determines the measure $\Prob$ uniquely. 
If, for a bounded Borel set $B\subset E$, we  let
$\#_B\colon\Conf(E)\to\mathbb{N}\cup\{0\}$ be the function
that to a configuration
assigns the number of its particles
belonging to~$B$, then, for any 
pairwise disjoint bounded Borel sets $B_1,\dotsc,B_l\subset E$
and any  $z_1,\dotsc,z_l\in {\mathbb C}$ from (\ref{eq1}) we  have
%\begin{equation}
%\label{eq2}
$\mathbb{E}_{\mathbb{P}}z_1^{\#_{B_1}}\dotsb z_l^{\#_{B_l}}
=\det\biggl(1+\sum\limits_{j=1}^l(z_j-1)\chi_{B_j}K\chi_{\sqcup_i B_i}\biggr).$
%\end{equation}

For further results and background on determinantal point processes, see e.g. \cite{BorOx},  \cite{HoughEtAl}, \cite{Lyons}, 
\cite{LyonsSteif}, \cite{Lytvynov}, \cite{ShirTaka0},  \cite{ShirTaka1}, \cite{ShirTaka2}, \cite{Soshnikov}.

In what follows we suppose that $K$ belongs to
$\scrI_{1, \text{loc}}(E,\mu)$, and denote
the corresponding determinantal measure by
$\mathbb{P}_K$. Note that $\mathbb{P}_K$ is uniquely defined by~$K$,
but different operators may yield the same measure.
By the Macch{\` \i}---Soshnikov theorem~\cite{Macchi}, \cite{Soshnikov}, any
Hermitian positive contraction that belongs
to the class~$\scrI_{1, \text{loc}}(E,\mu)$ defines a determinantal point process.

\subsection{Multiplicative functionals}

At the centre of the construction of infinite determinantal measures lies the result of \cite{Buf-umn} that can informally be summarized as follows:
a determinantal measure times a multiplicative functional is again a determinantal measure.
In other words, if $\mathbb{P}_K$ is a determinantal measure on $\Conf(E)$ induced by the
operator $K$ on $L_2(E,\mu)$, then, under certain additional assumptions,
it is shown in \cite{Buf-umn} that the measure ${ \Psi_g\mathbb{P}_K}$
after normalization yields a determinantal measure.

It is required in \cite{Buf-umn} that the operator $(g-1)K$ be of trace class; this assumption is too restrictive for our purposes, and in Propositions \ref{pr1-bis} and \ref{pr1} we shall now formulate two more convenient versions of Proposition 1 in \cite{Buf-umn}.

As before, let $g$ be a non-negative measurable function on $E$. If the operator $1+(g-1)K$ is invertible, then we set
$$B(g, K)=g K(1+{(g-1)}K)^{-1},\qquad
\tilde{B}(g, K)= {\sqrt{g}}K(1+{(g-1)}K)^{-1}{\sqrt{g}}.$$
By definition, $B(g,K),\tilde{B}(g,K)\in \scrI_{1,\text{loc}}(E,\mu)$ since $K\in \scrI_{1,\text{loc}}(E,\mu)$, and, if $K$ is self-adjoint, then
so is $\tilde{B}(g,K)$.

In the case when $K$ is self-adjoint,  the following proposition generalizes Proposition 1 in \cite{Buf-umn}.

\begin{proposition}
\label{pr1-bis}
Let $K\in \scrI_{1,\mathrm{loc}}(E,\mu)$ be a self-adjoint positive contraction, and let $\mathbb{P}_K$ be
the corresponding determinantal measure on $\Conf(E)$. Let
$g$ be a nonnegative bounded measurable
function on~$E$ such that
\begin{equation}
\label{gkint}
\sqrt{g-1}K\sqrt{g-1}\in\scrI_{1}(E,\mu)
\end{equation}
and that the operator
$1+{(g-1)}K$ is invertible. Then
\begin{enumerate}
\item we have $\Psi_g\in L_1(\Conf(E),\mathbb{P}_K)$ and
\begin{equation*}
\int\Psi_g\,d\mathbb{P}_K=\det\Bigl(1+\sqrt{g-1}K\sqrt{g-1}\Bigr)>0;
\end{equation*}
\item the operators
$B(g, K), \tilde{B}(g, K)$ induce on $\Conf(E)$ a determinantal measure
$\mathbb{P}_{B(g, K)}=\mathbb{P}_{\tilde{B}(g, K)}$ satisfying
\begin{equation}
\mathbb{P}_{B(g, K)}=\frac{\displaystyle \Psi_g\mathbb{P}_K}
{\displaystyle \int\limits_{\Conf(E)}\Psi_g\,d\mathbb{P}_K}.
\end{equation}
\end{enumerate}
\end{proposition}

{\bf Remark.} Since (\ref{gkint}) holds and $K$ is self-adjoint,  the operator $1+{(g-1)}K$ is 
invertible if and only if the operator  $1+\sqrt{g-1}K\sqrt{g-1}$ is invertible.

%By definition, the operator $\tilde B(g,K)$ is self-adjoint since so is $K$; also we have $\tilde B(g,K)\in\scrI_{1,\mathrm{loc}}(E,\mu)$ since $K\in\scrI_{1,loc}(E,\mu)$.

If $Q$ is a projection operator, then the operator $\tilde B(g,Q)$ admits the following description.

\begin{proposition}\label{projgl}
Let $L\subset L_2(E,\mu)$ be a closed subspace, and let
$Q$ be the operator of orthogonal projection onto $L$.
Let $g$ be a bounded measurable function such that the operator $1+(g-1)Q$ is invertible. Then the operator $\tilde B(g,Q)$ is the operator of orthogonal projection onto the closure of the subspace $\sqrt{g}L$.
\end{proposition}

We now consider the particular case when  $g$ is a characteristic function of a Borel subset.
In much the same way as before, if 
$E'\subset E$ is a Borel subset such that  the subspace $\chi_{E'}L$ is closed (recall that a sufficient condition for that is provided in Proposition \ref{closedsubsp}),  then we
set $Q^{E'}$ to be
the operator of orthogonal projection onto the closed subspace $\chi_{E'}L$.

Propositions \ref{pr1}, \ref{pr1-bis} now yield the following
\begin{corollary}\label{indsubset}
Let $Q\in\scrI_{1,\mathrm{loc}}(E,\mu)$ be the operator of orthogonal projection
onto a closed subspace $L\in L_2(E,\mu)$. Let $E'\subset E$ be a Borel subset such that
$\chi_{E'}Q\chi_{E'}\in\scrI_1(E,\mu)$.
Then
\begin{equation*}
\mathbb{P}_{Q}(\Conf(E, E'))=\det(1-\chi_{E\setminus E'}Q\chi_{E\setminus E'}).
\end{equation*}
Assume, additionally, that for any function $\varphi\in L$, the equality $\chi_{E'}\varphi=0$ implies $\varphi=0$.
Then  the subspace $\chi_{E'}L$ is closed, and we have
$$\mathbb{P}_{Q}(\Conf(E, E'))>0,  \  Q^{E'}\in\scrI_{1,\mathrm{loc}}(E,\mu),$$
and
\begin{equation*}
\frac{\mathbb{P}_{Q}|_{\Conf(E,E')}}{\mathbb{P}_{Q}(\Conf(E, E'))}=\mathbb{P}_{Q^{E'}}.
\end{equation*}
\end{corollary}

The induced measure of a determinantal measure onto the subset of configurations all whose particles lie in $E^{\prime}$ is thus again a determinantal measure. In the case of a discrete phase space, related induced processes were considered by Lyons \cite{Lyons} and by Borodin and Rains \cite{BorRains}.

Corollary \ref{indsubset} implies Theorem \ref{infdet-he}.

\subsection{The space $\scrI_\xi$}

To prove Proposition~\ref{pr1-bis}, we consider a slightly more general algebra of operators $K$ for which the trace  $\tr K$ and the Fredholm determinant $\det(1+K)$ can be defined and shown to have the usual properties. The  space
$\scrI_{\xi}(E,\mu)$
is a modification of  the space ${\mathcal L}_{1|2}(H)$ introduced
by Borodin, Okounkov and Olshanski \cite{BorOkOlsh} and used also by Olshanski in \cite{Olsh0}.
 We proceed to precise formulations.

Take a countable partition $\xi$ of our space $E$ into disjoint bounded measurable
sets $E_n$, $n\in {\mathbb N}$. 
Introduce the sets
\begin{equation}
\{\xi>n\}=\bigcup\limits_{k=n+1}^{\infty} E_k; \ \{\xi<n\}=\bigcup\limits_{k=1}^{n-1} E_k.
\end{equation}

Informally,  $\xi$ is considered as a random variable taking integer values.

The subspace   $$\scrI_{\xi}(E,\mu)\subset \scrI_{1, loc}(E,\mu)$$ is now defined as follows:
an operator $K\in \scrI_{1, loc}(E,\mu)$ belongs to  $\scrI_{\xi}(E,\mu)$ if
\begin{enumerate}
\item $K\in  \scrI_{2}(E,\mu)$;
\item $\sum\limits_{n=1}^{\infty} ||\chi_{E_n}K\chi_{E_n}||_{\scrI_{1}}<+\infty$.
\end{enumerate}
The space   $\scrI_{\xi}(E,\mu)$ is normed by the formula
$$
||K||_{\scrI_{\xi}}= ||K||_{\scrI_{2}}+\sum\limits_{n=1}^{\infty} ||\chi_{E_n}K\chi_{E_n}||_{\scrI_{1}}.
$$

By definition, the space  $\scrI_{\xi}(E,\mu)$ is an algebra. For $K\in \scrI_{\xi}(E,\mu)$, the Fredholm determinant
$\det(1+K)$ is defined by the formula
\begin{equation}
\label{def-det}
\det(1+K)=\det\left((1+K)\exp(-K)\right)\exp\left(\sum\limits_{n=1}^{\infty}
\tr(\chi_{E_n}K\chi_{E_n}) \right).
\end{equation}

The right-hand side of (\ref{def-det}) is well-defined since $(1+K)\exp(-K)\in {\scrI_{1}}$
for any $K\in {\scrI_{2}}$.

For $K_1, K_2\in {\scrI_{\xi}}$, we clearly have
$$
\det((1+K_1)(1+K_2))=\det(1+K_1)\det(1+K_2).
$$

From the definitions we now immediately obtain
\begin{proposition}
If
${(g-1)}K\in\scrI_{\xi}(E,\mu),$
then  $\Psi_g\in
L_1(\Conf(E),\mathbb{P}_K)$ and
$$\mathbb{E}_{\mathbb{P}_K}\Psi_g=\det(1+{(g-1)}K).$$
\end{proposition}

The following Proposition is a generalization of Proposition 1 in \cite{Buf-umn}.

\begin{proposition}
\label{pr1}
Assume that an operator $K\in
\scrI_{1,{loc}}(E,\mu)$ induces a
determinantal measure $\mathbb{P}_K$ on
$\Conf(E)$. Let $\xi$ be a countable measurable partition of $E$ and let
$g$ be a nonnegative bounded measurable
function on~$E$ such that ${(g-1)}K\in\scrI_{\xi}(E,\mu)$ and that the operator
$1+{(g-1)}K$ is invertible. Then the operators
$B(g, K), \tilde{B}(g, K)$ induce on $\Conf(E)$ a determinantal measure
$\mathbb{P}_{B(g, K)}=\mathbb{P}_{\tilde{B}(g, K)}$ satisfying
\begin{equation}
\displaystyle{\mathbb{P}_{B(g, K)}}=\frac{\displaystyle{\Psi_g\mathbb{P}_K}}
{\displaystyle{\int\limits_{\Conf(E)}}\displaystyle{\Psi_g d{\mathbb{P}_K}}}.
\end{equation}
\end{proposition}

Indeed, take a bounded measurable function $f$ on $E$ such that
$$(f-1)K\in\scrI_{\xi}(E,\mu).$$
We then immediately have
$$
\frac{{\mathbb E}_{{\mathbb P}_K} \Psi_f\Psi_g}{{\mathbb E}_{{\mathbb P}_K}\Psi_g}
=\det(1+(f-1)B(g,K))=
\det(1+{(f-1)}{\tilde B}(g,K)),
$$
and the proposition follows.

Observe now
that to a nonnegative function $g$ such that (\ref{gkint}) holds, one can easily assign a countable partition $\xi$ such that ${(g-1)}K\in\scrI_{\xi}(E,\mu)$. Proposition \ref{pr1-bis} is therefore clear
 from Proposition \ref{pr1}.

\section{Unitarily-Invariant  Measures on Spaces of Infinite Matrices}

\subsection{Pickrell Measures}

 Let $\Mat(n, \mathbb C)$ be the space of $n \times n$ matrices with complex entries:
\[
\Mat (n, \mathbb C) = \{ z = (z_{ij}),\; i=1,\ldots,\, n;j=1,\ldots, n \}
\]
Let $\mathrm{Leb}=dz$ be the Lebesgue measure on $\Mat (n, \mathbb C)$.

Following Pickrell \cite{Pickrell1}, take $s\in\mathbb{R}$ and introduce a measure $\widetilde{\mu}_n^{(s)}$ on $\mathrm{Mat}(n,\mathbb{C})$ by the formula $$\widetilde{\mu}_n^{(s)}=\det(1+{z}^*{z})^{-2n-s}dz.$$

The measure $\widetilde{\mu}_n^{(s)}$ is finite if and only  if $s>-1$.

For $n_1<n$, let $$\pi^n_{n_1}:\ \mathrm{Mat}(n,\mathbb{C})\to\mathrm{Mat}(n_1,\mathbb{C})$$ be the natural projection map that to a matrix ${z}=({z}_{ij}),i,j=1,\dots,n,$ assigns its upper left corner, the matrix $\pi_{n_1}^n({z})=({z}_{ij}),i,j=1,\dots,n_1.$

The measures $\widetilde{\mu}_n^{(s)}$ have the property of consistency with respect to the projections $\pi_{n_1}^n$.
More precisely, following Borodin and Olshanski \cite{BorOlsh}, p.116, observe that even if the measure $\widetilde{\mu}_n^{(s)}$ is
infinite, the fibres of the projection $\pi_{n-1}^n$ have finite conditional measure as long as $n+s>0$.
The push-forward $(\pi_{n-1}^n)_*\widetilde{\mu}_n^{(s)}$ is consequently well-defined, and
for any $s\in\mathbb{R}$ and $n>-s$ we have
\begin{equation} \label{projn}
(\pi_{n-1}^n)_*\widetilde{\mu}_n^{(s)}=\frac{\pi^{2n-1}(\Gamma(n+s))^2}{\Gamma(2n+s)\cdot\Gamma(2n-1+s)}\widetilde{\mu}_{n-1}^{(s)}.
\end{equation}

Now let $\mathrm{Mat}(\mathbb{N},\mathbb{C})$ be the space of infinite matrices whose rows and columns are indexed by natural numbers and whose entries are complex:
$$\mathrm{Mat}(\mathbb{N},\mathbb{C})=\{z=(z_{ij}),i,j\in\mathbb{N},z_{ij}\in\mathbb{C}\}.$$
Let $\pi_n^{\infty}:\mathrm{Mat}(\mathbb{N},\mathbb{C})\to\mathrm{Mat}(n,\mathbb{C})$ be the natural projection map that to an infinite matrix $z\in \mathrm{Mat}(\mathbb{N},\mathbb{C})$ assigns its  upper left $n\times n$-``corner'', the matrix $(z_{ij}), i,j=1,\dots,n.$

Take $s\in\mathbb{R}$ and $n_0\in\mathbb{N},n_0>-s.$ The relation (\ref{projn}) and the Kolmogorov Existence Theorem \cite{Kolmogorov} imply (for a detailed presentation, see  p. 116 in Borodin and Olshanski \cite{BorOlsh})  that for any $\lambda>0$ there exists a unique measure $\mu^{(s,\lambda)}$ on $\mathrm{Mat}(\mathbb{N},\mathbb{C})$ such that for any $n>n_0$ we have
\begin{equation}\label{infprojn}
(\pi_{n}^{\infty})_*\mu^{(s,\lambda)}=\lambda\left(\prod\limits_{l=n_0}^n\pi^{-2n}
\frac{\Gamma(2l+s)\Gamma(2l-1+s)}{(\Gamma(l+s))^2}\right)\widetilde{\mu}^{(s)}.
\end{equation}

If $s>-1$, the measures $\mu^{(s,\lambda)}$ are finite, and we let $\mu^{(s)}$ be the probability measure in the family $\mu^{(s,\lambda)}$.

In this case, (\ref{projn}) implies the relation
$$(\pi_{n}^{\infty})_*\mu^{(s)}=\pi^{-n^2}\prod\limits_{l=1}^n
\frac{\Gamma(2l+s)\Gamma(2l-1+s)}{(\Gamma(l+s))^2}\widetilde{\mu}_n^{(s)}.$$

If $s\leqslant-1$, the measures $\mu^{(s,\lambda)}$ are all infinite. In this case, slightly abusing notation, we shall omit the super-script $\lambda$ and write $\mu^{(s)}$ for a measure defined up to a multiplicative constant.
\begin{proposition} \label{singul}
For any $s_1, s_2\in {\mathbb R}$, $s_1\neq s_2$,  the Pickrell measures $\mu^{(s_1)}$ and 
$\mu^{(s_2)}$ are mutually singular.
\end{proposition}
Proposition \ref{singul}  is obtained from Kakutani's Theorem in the spirit of \cite{BorOlsh}, see also \cite{Neretin}.

Let $U(\infty)$ be the infinite unitary group: an infinite matrix $u=(u_{ij})_{i,j\in {\mathbb N}}$ belongs to
$U(\infty)$ if there exists a natural number $n_0$ such that the matrix
$$
(u_{ij}), i,j\in [1,n_0]
$$
is unitary, while $u_{ii}=1$ if $i>n_0$ and $u_{ij}=0$ if $i\neq j$, $\max(i,j)>n_0$.

The group $U(\infty) \times U(\infty)$ acts on $\Mat({\mathbb N}, {\mathbb C})$
by multiplication on both sides:
\[
T_{(u_1,u_2)}z \;=\; u_1zu_2^{-1}.
\]

The Pickrell measures $\mu^{(s)}$ are by definition $U(\infty)\times U(\infty)$-invariant.
For the r{\^o}le of Pickrell and related mesures in the representation theory of $U(\infty)$, see 
\cite{Olsh1}, \cite{Olsh2}, \cite{OlshVershik}.

Theorem 1 and Corollary 1 in \cite{Buf-ergdec} imply that the measures $\mu^{(s)}$ admit an ergodic decomposition. Furthermore,
Theorem 1 in \cite{Buf-inferg} implies that for any $s\in\mathbb{R}$ the ergodic components of the measure $\mu^{(s)}$ are almost surely finite. The main result of this note is an explicit description of the ergodic decomposition of the measures $\mu^{(s)}$ for $s\neq-1-2k,k\in\mathbb{N};$ in particular, for $s<-1$ we shall see that the ergodic decomposition is given by an explicitly computed infinite determinantal measure.

\subsection{Classification of ergodic measures}

First, we recall the classification of ergodic probability $U(\infty)\times U(\infty)$-invariant measures on
$\Mat(\mathbb{N},\mathbb{C})$. This classification has been obtained
by Pickrell \cite{Pickrell1}, \cite{Pickrell2}; Vershik \cite{Vershik} and Olshanski and Vershik  \cite{OlshVershik} proposed a different approach to this classification in the case of unitarily-invariant measures on the space of infinite Hermitian matrices, and Rabaoui \cite{Rabaoui1}, \cite{Rabaoui2} adapted the Olshanski-Vershik approach to the initial problem of Pickrell.
In this note, the Olshanski-Vershik approach is followed as well.

Take $z\in \Mat(\mathbb{N},\mathbb{C})$, denote $z^{(n)}=\pi^{\infty}_{n}z,$
and let 
\begin{equation}\label{lan}
\lambda_1^{(n)}\geqslant\ldots\geqslant\lambda_n^{(n)}\geqslant0
\end{equation} 
be the eigenvalues of
the matrix $$\left(z^{(n)}\right)^*z^{(n)},$$
counted with multiplicities, arranged in non-increasing order.
To stress dependence on $z$, we write $\lambda_i^{(n)}=\lambda_i^{(n)}(z)$.

\begin{teorema}
\begin{enumerate}
\item Let $\eta$ be an ergodic Borel $U(\infty)\times U(\infty)$-invariant probability measure on $\Mat(\mathbb{N},\mathbb{C})$. Then there exist non-negative real numbers
$$\gamma\geqslant0,\;\,x_1\geqslant x_2\geqslant\ldots\geqslant x_n\geqslant\ldots\geqslant0\,,$$
satisfying $\displaystyle\gamma\geqslant\sum\limits_{i=1}^{\infty}x_i$, such that for $\eta$-almost every
$z\in\mathrm{Mat}(\mathbb{N},\mathbb{C})$ and any $i\in\mathbb{N}$ we have:
\begin{equation}\label{xgamma}
x_i=\lim\limits_{n\to\infty}\frac{\lambda_i^{(n)}(z)}{n^2},\ \
\gamma=\lim\limits_{n\to\infty}\frac{\tr\left(z^{(n)}\right)^*z^{(n)}}{n^2}.
\end{equation}
\item Conversely, given non-negative real numbers
$\gamma\geqslant0,\;\,x_1\geqslant x_2\geqslant\ldots\geqslant x_n\geqslant\ldots\geqslant0$ such that
$$\displaystyle\gamma\geqslant\sum\limits_{i=1}^{\infty}x_i\,,$$
there exists a unique $U(\infty)\times U(\infty)$-invariant ergodic Borel probability measure $\eta$ on
$\mathrm{Mat}(\mathbb{N},\mathbb{C})$ such that the relations (\ref{xgamma}) hold for $\eta$-almost all $z\in \Mat(\mathbb{N},\mathbb{C})$.
\end{enumerate}
\end{teorema}

Introduce the {\it Pickrell set} $\Omega_P\subset\mathbb{R}_+\times\mathbb{R}_+^\mathbb{N}$ by the formula
$$\Omega_P=\left\{\omega=(\gamma,x)\colon x=(x_n),\;n\in\mathbb{N},\;x_n\geqslant x_{n+1}\geqslant0,\;
\gamma\geqslant\sum\limits_{i=1}^{\infty}x_i\right\}.$$
The set $\Omega_P$ is, by definition, a closed subset of $\mathbb{R}_+\times\mathbb{R}_+^\mathbb{N}$ endowed with the Tychonoff topology.

By Proposition 3 in \cite{Buf-ergdec}, the subset of ergodic $U(\infty)\times U(\infty)$-invariant measures
is a Borel subset of the space of all Borel probability measures on $\mathrm{Mat}(\mathbb{N},\mathbb{C})$ endowed
with the natural Borel structure (see, e.g., \cite{Bogachev}).
Furthermore, if one denotes $\eta_{\omega}$ the Borel ergodic probability measure corresponding to a point $\omega\in\Omega_P$, $\omega=(\gamma,x)$, then the correspondence
$$\omega\longrightarrow\eta_{\omega}$$
is a Borel isomorphism of the Pickrell set $\Omega_P$ and the set of $U(\infty)\times U(\infty)$-invariant ergodic probability measures on $\mathrm{Mat}(\mathbb{N},\mathbb{C})$.

The Ergodic Decomposition Theorem (Theorem 1 and Corollary 1 of \cite{Buf-ergdec})  implies that
each Pickrell measure $\mu^{(s)}$, $s\in\mathbb{R}$, induces a unique decomposing measure $\overline{\mu}^{(s)}$ on $\Omega_P$ such that we have
\begin{equation}
\label{ergdecmubar}
\mu^{(s)}=\int\limits_{\Omega_P}\eta_{\omega}\,d\overline{\mu}^{(s)}(\omega)\;.
\end{equation}
The integral is understood in the usual weak sense, see \cite{Buf-ergdec}.

For $s>-1$, the measure $\overline{\mu}^{(s)}$ is a probability measure on $\Omega_P$,
while for $s\leqslant-1$ the measure $\overline{\mu}^{(s)}$ is infinite.

Set
$$\Omega_P^0=\{ (\gamma,\{x_n\})\in\Omega_P: x_n>0\text{\quad for all }n,\  \gamma=\sum\limits_{n=1}^{\infty} x_n \}.$$

The subset $\Omega_P^0$ is of course { not} closed in $\Omega_P$.

Introduce a map $$
\conf\colon\Omega_P\to  \Conf((0, +\infty))
$$
that to a point $\omega\in\Omega_P, \omega=(\gamma, \{x_n\})$ assigns
the configuration $$
\conf(\omega)=(x_1, \ldots, x_n, \ldots)\in
\Conf((0, +\infty)).
$$
The map $\omega\to \conf(\omega)$ is bijective in restriction to the subset
$\Omega_P^0$.

{\bf Remark.} In the definition of the map $\conf$, the ``asymptotic eigenvalues'' $x_n$ are  counted with multiplicities, while, if $x_{n_0}=0$ for some $n_0$, then $x_{n_0}$
and all subsequent terms are discarded, and the resulting configuration is finite. We shall see, however,
that the complement $\Omega_P\backslash\Omega_P^0$ is $\overline{\mu}^{(s)}$-negligible for all $s\neq -1-2k$, $k\in {\mathbb N}$,  and, consequently,
that, $\overline{\mu}^{(s)}$-almost surely, all configurations are infinite. It will also develop that,
$\overline{\mu}^{(s)}$-almost surely, all multiplicities are equal to one.

We proceed to the formulation of the  main result of this note,  an explicit description of the measures $\overline{\mu}^{(s)}$ for
$s\neq -1-2k$, $k\in{\mathbb N}$.

\section{  Ergodic decomposition of infinite Pickrell measures}

\subsection{The Bessel point process and finite Pickrell measures.}

Consider the half-line $(0, +\infty)$ endowed with the standard Lebesgue measure $\mathrm{Leb}$.
Take $s>-1$ and consider the standard Bessel kernel
\begin{equation}\label{bessel}
J_s(x,y)=\frac{\sqrt{x}J_{s+1}(\sqrt{x})J_s(\sqrt{y})-\sqrt{y}J_{s+1}(\sqrt{y})J_s(\sqrt{x})}{2(x-y)}
\end{equation}
(see, e.g., page 295 in Tracy and Widom \cite{TracyWidom}). The kernel $J_s$ induces on $L_2((0, +\infty), {\mathrm{Leb}})$ the operator of orthogonal projection onto the subspace of functions whose
 Hankel transform is supported in $[0,1]$ (see \cite{TracyWidom}). 
Setting 
$
x_1=4/x, \ x_2=4/y
$
yields a kernel $K^{(s)}$ given by the formula
\begin{multline}
K^{(s)}(x_1,x_2)=\frac{J_s\left(\frac{2}{\sqrt{x_1}}\right)\frac{1}{\sqrt{x_2}}J_{s+1}\left(\frac{2}{\sqrt{x_2}}\right)- J_s\left(\frac{2}{\sqrt{x_2}}\right)\frac{1}{\sqrt{x_1}}J_{s+1}\left(\frac{2}{\sqrt{x_1}}\right)}{x_1-x_2}, \\ x_1>0,\;x_2>0\,.
\end{multline}
(recall here that a change of variables $u_1=\rho(v_1)$, $u_2=\rho(v_2)$ transforms a kernel $K(u_1, u_2)$
to a kernel of the form $K(\rho(v_1), \rho(v_2))(\sqrt{\rho^{\prime}(v_1)\rho^{\prime}(v_2)})$).

The kernel $K^{(s)}$ induces on the space $L_2((0,+\infty), \mathrm{Leb})$
a locally trace class operator of orthogonal projection, for which, slightly abusing notation, we keep the symbol $K^{(s)}$;
 by the Macch\`{\i}-Soshnikov Theorem, the operator $K^{(s)}$
 induces a determinantal measure $\mathbb{P}_{K^{(s)}}$ on $\Conf((0,+\infty))$.
The determinantal measure $\mathbb{P}_{K^{(s)}}$ is precisely the decomposing measure for the Pickrell measure $\mu^{(s)}$, as is shown
by the following
\begin{proposition}\label{finpick} Let $s>-1$. Then $\overline{\mu}^{(s)}(\Omega_P^0)=1$ and
the  $\overline{\mu}^{(s)}$-almost sure bijection
$\omega\to \conf(\omega)$ identifies the measure $\overline{\mu}^{(s)}$ with the determinantal measure $\mathbb{P}_{K^{(s)}}$.
\end{proposition}
Sketch of proof of Proposition \ref{finpick}.
Take $s>-1$ and let $P_n^{(s)}(u)$ be the standard Jacobi orthogonal
polynomials on the interval $[-1,1]$,
corresponding to the weight $(1-u)^s$.

Following Pickrell,  to a matrix $z\in Mat(n,\mathbb{C})$ assign the collection
$\left(\lambda_1(z),\ldots,\lambda_n(z)\right)$ of the
eigenvalues of the matrix $z^*z$ arranged in non-increasing order  ( cf. (\ref{lan})). 
The {\it radial part} $\mathfrak{r}^{(n,s)}$ 
of the Pickrell 
measure $\mu_n^{(s)}$ is now defined as the 
push-forward of the  measure
$\mu_n^{(s)}$ under the
map $$z\to\left(\lambda_1(z),\ldots,\lambda_n(z)\right).$$ The radial part of the Pickrell measure has
determinantal form:
$$d\mathfrak{r}^{(n,s)}(\lambda)=\frac1{n!}\det
K_n^{(s)}\left(\lambda_i,\lambda_j\right)\;\prod\limits_{i=1}^n\,d\lambda_i,\;\;\lambda_i>0\,.$$
where
\begin{multline}
K_n^{(s)}\left(\lambda_1, \lambda_2\right)=\frac{n(n+s)}{(2n+s)(1+\lambda_1)^{s/2}(1+\lambda_2)^{s/2}} \times\\
\times \frac{
P_n^{(s)} \left( \frac{\lambda_1-1}{\lambda_1+1} \right)
P_{n-1}^{(s)} \left( \frac{\lambda_2-1}{\lambda_2+1} \right) - 
P_n^{(s)} \left( \frac{\lambda_2-1}{\lambda_2+1} \right)
P_{n-1}^{(s)} \left( \frac{\lambda_1-1}{\lambda_1+1} \right)
}{\lambda_1-\lambda_2}.
\end{multline}
The change of variables $\displaystyle
u_i=\frac{\lambda_i-1}{\lambda_i+1}$, $i=1, \dots, n$,  reduces $K_n^{(s)}$ to the
Christoffel-Darboux kernel for the Jacobi orthogonal ensemble with weight
$(1-u)^s$.

Introducing the scaling $\lambda_i=n^2x_i$, taking $n\to\infty$ and using
the classical asymptotics for Jacobi orthogonal
polynomials (see, e.g., Szeg{\"o}  \cite{Szego}), one finds
$$\lim\limits_{n\to\infty}n^2K_n^{(s)}\left(n^2x_1,n^2x_2\right)=K^{(s)}(x_1,x_2)\,,$$
convergence being uniform on compact subsets of $(0,+\infty)$. 
To prove that $\overline{\mu}^{(s)}(\Omega_P^0)=1$, the method of Section 7 in Borodin and Olshanski \cite{BorOlsh} is adapted to our situation.

\subsection{A recurrence relation for Bessel point proceses.} 
The following  observation motivates the construction of the next subsection.
Given a finite family of functions $f_1, \dots, f_N$ on the real line, let $\spann(f_1, \dots, f_N)$ stand for the vector space these functions span.
For any $s\in {\mathbb R}$, $N\in {\mathbb N}$ we clearly have 
\begin{multline}\label{rec-jac}
\spann\left((1-u)^{s/2}
, \dots, (1-u)^{s/2}u^N\right)={\mathbb R}(1-u)^{s/2}\oplus\\
\oplus \spann\left((1-u)^{(s+2)/2}
, \dots, (1-u)^{(s+2)/2}u^{N-1}\right).
\end{multline}
If $s>-1$, then the informal meaning of (\ref{rec-jac}) is that the space of the first $N+1$ normalized Jacobi polynomials with weight $(1-u)^s$ is a rank one perturbation of the space of the first $N$  normalized Jacobi polynomials with weight $(1-u)^{s+2}$.

A similar statement holds true for the Bessel kernel: using the recurrence relation
%\begin{equation}\label{rec-bes}
$J_{s+1}(x)=\frac{2s}{x}J_s(x)-J_{s-1}(x)$
%\end{equation}
for Bessel functions,
one easily obtains the recurrence  relation
\begin{equation}
J_s(x,y)=J_{s+2}(x,y)+\frac{s+1}{\sqrt{xy}}J_{s+1}(\sqrt{x})J_{s+1}(\sqrt{y})
\end{equation}
for the Bessel kernels: the Bessel kernel with parameter $s$ is thus a rank one perturbation of the Bessel 
kernel with parameter $s+2$.

For ergodic decomposition measures of infinite Pickrell measures we shall now give a similar description 
in terms of infinite determinantal measures obtained as finite-rank perturbations of Bessel point processes.

\subsection{Formulation of the main result.}

Now take $s<-1$, $s\neq -1-2k$, $k\in {\mathbb N}$.
Let $n_s$ be such that
$$
\frac{s}{2}+n_s\in\left(-\frac12, \frac12\right).
$$

Introduce a finite-dimensional subspace $V^{(s)}\subset L_{2,loc}((0,+\infty), \mathrm{Leb})$ by the formula
$$V^{(s)}=\operatorname{span}\left(x^{-s\!/\!_2-1}, \ldots, x^{-s\!/\!_2-n_s}\right)\,.$$

For $s^{\prime}>-1$, let $L^{(s^{\prime})}\subset L_2((0, +\infty), \mathrm {Leb})$ be the range of the operator $K^{s^{\prime}}$, and for $s<-1$, $s\neq -1-2k$, $k\in {\mathbb N}$, introduce a subspace  $H^{(s)}$ of $L_{2, \mathrm{loc}}((0, +\infty), \mathrm{Leb})$
by the formula
$$
H^{(s)}=L^{(s+2n_s)}+V^{(s)}.
$$
%Given $s<-1$, $s\neq -1-2k$, $k\in {\mathbb N}$, u
Using Proposition \ref{lvehe}, one easily checks that if $R>0$ is big enough, then the subspace $H^{(s)}\subset L_2((0, +\infty), \mathrm {Leb})$ and the subset $E_0=(0,R)$ satisfy Assumption \ref{he}. Let
\begin{equation}\label{bs-def}
%$$
{\mathbb B}^{(s)}={\mathbb B}(H^{(s)}, (0,R))
 %$$
 \end{equation}
 be the corresponding infinite determinantal measure (which, by definition, does not depend on the specific choice of a big enough $R$).

The ergodic decomposition of infinite Pickrell measures is now given by the following

\begin{theorem} \label{mainthm}
Let $s<-1$, $s\neq -1-2k$, $k\in {\mathbb N}$, and let $\overline{\mu}^{(s)}$
be the decomposing measure, defined by (\ref{ergdecmubar}), of the Pickrell measure ${\mu}^{(s)}$.
Then \begin{enumerate}
\item $\overline{\mu}^{(s)}(\Omega_P\backslash \Omega_P^0)=0$;

\item the $\overline{\mu}^{(s)}$-almost sure bijection $\omega\to {\rm conf}(\omega)$ identifies
$\overline{\mu}^{(s)}$ with the infinite determinantal measure
$\mathbb{B}^{(s)}$.
%defined by (\ref{bs-def}).
\end{enumerate}
\end{theorem}
Take $R>0$ and set $$\Omega_P^0(R)=\{\omega\in\Omega_P^0: x_1\leq R\}.$$ 
Set $L_R^{(s)}=\chi_{(0,R)}H^{(s)}$; the subspace  $L_R^{(s)}$ is closed if $R>0$ is big enough,
and we let $Q_R^{(s)}$ be the corresponding operator of orthogonal projection.
By Proposition \ref{gn-conv}, we have  $Q_R^{(s)}\to K^{(s+2n_s)}$ in ${\mathscr I}_{1, {\mathrm loc}}((0, +\infty), {\mathrm {Leb}})$
as $R\to\infty$.
Theorem \ref{mainthm} together with Theorem \ref{infdet-he} implies 
%$\overline{\mu}^{(s)}(\Omega_P^0(R))<+\infty$
\begin{corollary} If $R$ is big enough, then  
$\overline{\mu}^{(s)}(\Omega_P^0(R))<+\infty$
and the $\overline{\mu}^{(s)}$-almost sure bijection $\omega\to {\rm conf}(\omega)$ identifies the normalized restriction of the measure
$\overline{\mu}^{(s)}$ to the subset ${\Omega_P^0(R)}$ with the determinantal measure
$\mathbb{P}_{Q_R^{(s)}}$.
\end{corollary}
The proof of Theorem \ref{mainthm} starts, again, with the computation of the radial part of the infinite Pickrell measure; changing variables by the formula 
$$
u_i=\frac{\lambda_i-1}{\lambda_i+1}, \ i=1, \dots, n,
$$ 
one arrives at an ``infinite orthogonal polynomial ensemble''(cf. (\ref{ore})) of the form
\begin{equation}\label{inforth}
\prod\limits_{i<j}(u_i-u_j)^2\prod\limits_{i} (1-u_i)^s.
\end{equation}
By definition, the measure (\ref{inforth}) is an infinite determinantal measure 
obtained by perturbing the closed subspace
$$
\spann\left((1-u)^{(s+2n_s)/2}, \dots, (1-u)^{(s+2n_s)/2}u^{N-n_s-1}\right)\subset L_2([0,1], {\mathrm {Leb}})
$$
by the finite-dimensional subspace 
$$
\spann\left( (1-u)^{s/2}, \dots,         (1-u)^{(s+2n_s-2)/2}\right)\subset L_{2, {\mathrm loc}}([0,1], {\mathrm {Leb}}).
$$
The next step is to take the scaling limit of these infinite determinantal measures. This is achieved, with the use of Propositions \ref{mult-he} and \ref{finite-mult},  by taking the product with a suitably chosen multiplicative functional and effecting the scaling limit transition for the corresponding determinantal probability measures.
The detailed proof of Theorem \ref{mainthm} will be published in the sequel to this note.

{\bf{Acknowledgements.}}
Grigori Olshanski posed the problem to me, and I am greatly indebted to him.
I am deeply grateful to Alexei Borodin, Alexei Klimenko, Klaus Schmidt  and Alexander Soshnikov for
many useful discussions. I am deeply grateful to the referees of the paper for many helpful suggestions.
I am deeply grateful to Sergei Sharahov and Ramil Ulmaskulov for typesetting parts of the manuscript.

This work has been supported in part by an Alfred P. Sloan Research Fellowship, a
Dynasty Foundation Fellowship, as well as an IUM-Simons Fellowship,
by the Grant MK-6734.2012.1 of the President of the Russian Federation,
by the Programme ``Dynamical systems and mathematical control theory''
of the Presidium of the Russian Academy of Sciences,
by the RFBR-CNRS grant 10-01-93115-NTsNIL and by the RFBR grant 11-01-00654.

During work on this paper I was visiting the Institut Mittag-Leffler of
the Royal Swedish Academy of Sciences in Djursholm, the Institut Henri Poincar{\'e} in Paris,
the Centre International de Rencontres Math{\'e}matiques  in Marseille, the Abdus Salam International Centre for
Theoretical Physics in Trieste,  the Joint Institute for Nuclear Research in Dubna, the University of Kyoto and the
Chebyshev Laboratory in Saint-Petersburg.
I am deeply grateful to these institutions for their warm hospitality.

\end{document}